\documentclass[a4paper,10pt,,reqno]{amsart}
\usepackage[utf8]{inputenc}
\usepackage{amssymb}

\usepackage{color}

\usepackage[colorlinks=true]{hyperref}

\definecolor{dark-red}{rgb}{.54,.0,.0}
\definecolor{dark-green}{rgb}{.0,.4,.0}
\definecolor{dark-blue}{rgb}{.04,.04,.4}

\hypersetup{linkcolor=dark-red, urlcolor=dark-blue, citecolor=dark-green}

\renewcommand{\theequation}{\arabic{equation}}

\renewcommand{\wp}{\mbox{$\omega_{i+1}$}}
\renewcommand{\d}{\displaystyle}
\newcommand{\cSection}{Section}
\newcommand{\cSectiontwo}{Section~2}
\newcommand{\cpaper}{paper}
\newcommand{\n}{\newcommand}
\newcounter{transfer}
\n{\m}{\mbox}
\n{\w}{\m{$\omega$}}
\n{\W}{\m{$\varpi$}}
\n{\equal}{\m{$\stackrel{\triangle}{=}$}}
\n{\fpp}{\m{$f^{\prime\prime}$}}
\n{\ffi}{\m{$f_i$}}
\n{\ffm}{\m{$f_{i-1}$}}
\n{\ffp}{\m{$f_{i+1}$}}
\n{\dt}{\m{$\Delta\theta$}}
\n{\cosdt}{\m{$\cos\dt$}}
\n{\tomc}{\m{$2 (1 - \cosdt)$}}
\n{\fsecdrv}{\m{$\d\frac{\ffp - 2 \ffi + \ffm}{\tomc}$}}
\n{\wi}{\m{$\w_i$}}
\n{\wm}{\m{$\w_{i-1}$}}
\n{\wid}{\m{$\dot{\w}_i$}}
\n{\wsecdrv}{\m{$\d\frac{\wp - 2 \wi + \wm}{\tomc}$}}
\n{\wic}{\m{$\wi^3$}}
\newcommand{\change}{\newcommand}

\newcounter{initialplyg}

\renewcommand{\l}{\mbox{$\log\,$}}

\change{\half}{\m{$\d\frac{1}{2}$}}
\change{\ddt}{\m{$\d\frac{d}{dt}$}}
\n{\be}{\begin{equation}}
\n{\ee}{\end{equation}}
\n{\pp}{\m{$p$}}
\n{\tpi}{\m{$2\pi$}}
\n{\dtau}{\m{$d\tau$}}
\n{\dvt}{\m{$d\vartheta$}}
\n{\odt}{\m{$\m{O}(\dt)$}}
\n{\dts}{\m{$(\dt)^2$}}
\n{\odts}{\m{$\m{O}\!\left [(\dt)^2\right ]$}}
\n{\odtc}{\m{$\m{O}\!\left [(\dt)^3\right ]$}}
\n{\odtf}{\m{$\m{O}\!\left [(\dt)^4\right ]$}}
\n{\tmidt}{\m{$(\theta - i\dt)$}}
\n{\otmidts}{\m{$\m{O}\!\left [(\theta - i\dt)^2\right ]$}}
\n{\ovtmidts}{\m{$\m{O}\!\left [(\vartheta - i\dt)^2\right ]$}}
\n{\itm}{\m{$\d\frac{1}{2M}$}}
\n{\sindt}{\m{$\sin\dt$}}
\n{\tandt}{\m{$\tan\dt$}}
\n{\cotdt}{\m{$\cot\dt$}}
\n{\cscdt}{\m{$\csc\dt$}}
\n{\omcpm}{\m{$1 - \cos\frac{\pi}{M}$}}
\n{\domcpm}{\m{$\d 1 - \cos\frac{\pi}{M}$}}
\n{\fraction}{\m{$\d\frac{\tomc}{\sindt}$}}
\n{\ifraction}{\m{$\d\frac{\sindt}{\tomc}$}}
\n{\itomc}{\m{$\d\frac{1}{\tomc}$}}
\n{\almostone}{\m{$\d\frac{\sin^2\dt}{\tomc}$}}
\n{\closeone}{\m{$\d\frac{\sqrt{\tomc}}{\dt}$}}
\n{\bigli}{\m{$L_i$}}
\n{\ki}{\m{$k_i$}}
\n{\km}{\m{$k_{i-1}$}}
\n{\kp}{\m{$k_{i+1}$}}
\n{\gi}{\m{$g_i$}}
\n{\gm}{\m{$g_{i-1}$}}
\n{\gp}{\m{$g_{i+1}$}}
\n{\hi}{\m{$h_i$}}
\n{\hm}{\m{$h_{i-1}$}}
\n{\hp}{\m{$h_{i+1}$}}
\n{\kid}{\m{$\dot{k}_i$}}
\n{\kpi}{\m{$\kappa_i$}}
\n{\kpm}{\m{$\kappa_{i-1}$}}
\n{\kpp}{\m{$\kappa_{i+1}$}}
\n{\kis}{\m{$\ki^2$}}
\n{\kic}{\m{$\ki^3$}}
\n{\ikpi}{\m{$\d\frac{1}{\kpi(0)}$}}
\n{\ks}{\m{$k_\ast$}}
\n{\lidot}{\m{$2 \wi \cotdt - \wm \cscdt - \wp \cscdt$}}
\change{\Wi}{\m{$\W_i$}}
\change{\Wm}{\m{$\W_{i-1}$}}
\change{\Wp}{\m{$\W_{i+1}$}}
\n{\wis}{\m{$\wi^2$}}
\n{\iWi}{\m{$\d\frac{1}{\Wi(0)}$}}
\n{\ws}{\m{$\w_\ast$}}
\change{\li}{\m{$\Lambda_i$}}
\change{\lm}{\m{$\Lambda_{i-1}$}}
\change{\lp}{\m{$\Lambda_{i+1}$}}
\n{\upi}{\m{$\Upsilon_i$}}
\n{\upm}{\m{$\Upsilon_{i-1}$}}
\change{\upp}{\m{$\Upsilon_{i+1}$}}
\change{\ui}{\m{$u_i$}}
\change{\um}{\m{$u_{i-1}$}}
\change{\up}{\m{$u_{i+1}$}}
\n{\vi}{\m{$v_i$}}
\n{\vm}{\m{$v_{i-1}$}}
\n{\vp}{\m{$v_{i+1}$}}
\n{\ri}{\m{$r_i$}}
\n{\rp}{\m{$r_{i+1}$}}
\n{\pn}{\m{$P_n$}}
\n{\ppm}{\m{$p_m$}}
\n{\pmp}{\m{$p_{m+1}$}}
\n{\mj}{\m{$M_j$}}
\n{\ijnew}{\m{$i_j$}}
\n{\secdrv}{\m{$\d\frac{\kp - 2 \ki + \km}{\tomc}$}}
\n{\Wsecdrv}{\m{$\d\frac{\Wp - 2 \Wi + \Wm}{\tomc}$}}
\n{\secdrvkp}{\m{$\d\frac{\kpp - 2 \kpi + \kpm}{\tomc}$}}
\n{\secdrvu}{\m{$\d\frac{\up - 2 \ui + \um}{\tomc}$}}
\n{\secdrvv}{\m{$\d\frac{\vp - 2 \vi + \vm}{\tomc}$}}
\n{\secdrvl}{\m{$\d\frac{\lp - 2 \li + \lm}{\tomc}$}}
\n{\secdrvup}{\m{$\d\frac{\upp - 2 \upi + \upm}{\tomc}$}}
\n{\fdrv}{\m{$\d\frac{\kp - \ki}{\sindt}$}}
\n{\fdrvs}{\m{$\d\frac{(\kp - \ki)^2}{\tomc}$}}
\n{\fdrvc}{\m{$\d\frac{\kp - \km}{2 \sindt}$}}
\n{\fdrvcz}{\m{$\d\frac{k_1 - k_{-1}}{2 \sindt}$}}
\n{\threeks}{\m{$\ki^2 + \ki\kp + \kp^2$}}
\n{\twok}{\m{$\ki + \kp$}}
\n{\kmk}{\m{$\kpi - \ki$}}
\n{\wfdrv}{\m{$\d\frac{\wp - \wi}{\sindt}$}}
\n{\wfdrvs}{\m{$\d\frac{(\wp - \wi)^2}{\tomc}$}}
\n{\wfdrvc}{\m{$\d\frac{\wp - \wm}{2 \sindt}$}}
\n{\wfdrvcz}{\m{$\d\frac{w_1 - w_{-1}}{2 \sindt}$}}
\n{\threews}{\m{$\wi^2 + \wi\wp + \wp^2$}}
\n{\twow}{\m{$\wi + \wp$}}
\n{\wmw}{\m{$\Wi - \wi$}}
\n{\ep}{\m{$e^{2\pi imn/(2M)}$}} 
\n{\emmn}{\m{$e^{-2\pi imn/(2M)}$}} 
\n{\emn}{\m{$e^{-2\pi in/(2M)}$}} 
\n{\epp}{\m{$e^{2\pi iq/(2M)}$}} 
\n{\rav}{\m{$\d\frac{\rp + \ri}{2}$}}
\n{\ravz}{\m{$\d\frac{r_1 + r_0}{2}$}}
\n{\kptidt}{\m{$\kappa_{\theta}(i\dt,0)$}}
\n{\Wtidt}{\m{$\W_{\theta}(i\dt,0)$}}
\n{\optr}{\m{$1 + \tau \rho_i + \m{O}(\tau^2)$}}
\change{\ut}{\m{$u_\theta$}}
\n{\utt}{\m{$u_{\theta\theta}$}}
\change{\s}{\m{$\d\sum_i\ $}}
\n{\stm}{\m{$\d\sum_{n=1}^{2M-1}\ $}}
\n{\smtm}{\m{$\d\sum_{m=0}^{2M-1}\ $}}
\change{\sm}{\m{$\d\sum_{m=0}^{M-1}$}}
\n{\intzdt}{\m{$\d\int_0^{\Delta\theta}$}}
\n{\lip}{\m{$L_i^+$}}
\n{\limn}{\m{$L_i^-$}}
\n{\firstint}{\m{$\d\intzdt \frac{\cos\tau}{\optr} d\tau$}}
\n{\secint}{\m{$\d\intzdt \frac{\sin\tau}{\optr} d\tau$}}
\n{\firstintm}{\m{$\d\frac{1}{\kpi} \int_0^{\Delta\theta/2}
                \frac{\cos\tau}{\optr} d\tau$}}
\n{\secintm}{\m{$\d\frac{1}{\kpi} \int_0^{\Delta\theta/2}
                \frac{\sin\tau}{\optr} d\tau$}}

\title[Crystalline Algorithm for motion by curvature]{Convergence of 
a crystalline algorithm for the motion
of a simple closed convex curve by weighted curvature}

\author{Pedro Martins Gir\~{a}o}

\address{Courant Institute,
251 Mercer Street, New York, NY 10012}
\curraddr{Mathematics Department,
Instituto Superior T\'{e}cnico,
1049-001 Lisbon, Portugal}

\email{pgirao@math.ist.utl.pt}

\thanks{This work was partially supported by AFOSR grant 90-0090.}

\begin{document}

\begin{abstract}
Motion by weighted mean curvature is a geometric evolution law
for surfaces 
and represents
steepest descent with respect to
anisotropic surface energy. It has been proposed that this
motion could be computed 
numerically by using a ``crystalline'' approximation to the
surface energy in the evolution law.
In this \cpaper\ we prove the convergence of this
numerical method for the case 
of simple closed convex curves in the plane.
\end{abstract}

\keywords{crystalline, motion by curvature, surface-energy}

\subjclass{65M12, 73B30}

\maketitle

\section{Introduction}

In the modeling of phase 
transitions it is often of interest
to consider surface-energy-driven motion of interfaces.  If the
surface energy is anisotropic this leads to motion by 
weighted mean curvature
(see Angenent and Gurtin \cite{AG}).
Some materials have ``crystalline'' energies.
The idea of taking a crystalline 
approximation of a convex surface-energy leads naturally to a numerical 
scheme for determining motion by weighted curvature.
For a detailed explanation as well as for the physical and
mathematical context of this work the reader is referred to
the introduction and appendix of \cite{GK}.
(Familiarity with the theory of surface-energy-driven motion
of phase boundaries is {\em not} assumed either in this work
or in \cite{GK}). 
The idea of taking a crystalline approximation to the surface-energy
has been studied in the static case by Sullivan in \cite{S} to
find an area minimizing oriented hypersurface that spans a given
boundary in space.
The first attempt to study this idea in the dynamic context
was \cite{GK}, where we examined the case in which the interface is
the graph of a function of one variable.
The main result of \cite{GK} was the convergence of a certain numerical 
scheme for a quasi\-linear parabolic differential equation with 
constant Dirichlet or Neumann boundary conditions.  
We proved convergence in $H^1$\/ with a specified rate.
Our method was somewhat similar to the convergence analysis
for a Galerkin approximation.

Here we study the crystalline approximation of motion
by weighted curvature for smooth simple closed convex curves.
For motion by weighted curvature the normal velocity of
the curve is its weighted curvature; 
the weighted curvature $\W$\/ is
\[
        \W\ \equal\ (f + \fpp)\,\kappa ,
\]
if $f$\/ (the interfacial energy) is a smooth function which depends
only on the angle between the normal to the curve and a fixed axis
and $\kappa$\/ is the curvature.
We assume throughout that both $f$\/ and $f + \fpp$\/ are positive functions.
We approximate the interface by a polygon with $N$ ($\geq\ 4$) sides
moving by weighted curvature (as defined below).  
All the interior angles of the polygon are equal and
the sides of the polygon have fixed directions:
the interior normal to the $i$th side is
$N_i\ =\ -(\cos i\Delta\theta, \sin i\Delta\theta)$,
where $\Delta\theta\ =\ 2\pi/N$,
in a fixed coordinate system.

It is standard to  parametrize a convex curve by the angle 
$\theta$\/ between its
normal and a fixed coordinate axis (here $\theta$\/ is the
angle between the exterior normal and the $x$-axis).  
Using this parametrization
the evolution equation for the weighted curvature of the interface is
\[
        \W_t\ =\ (f + \fpp)^{-1} (\W^2 \W_{\theta\theta} + \W^3) 
\]
(see Eq.~(2.23) of Angenent and Gurtin~\cite{AG}).
On the other hand,
the evolution equation for the weighted curvature 
$\w_i$\/ of the $i$th side of the polygon turns out to be
\[
        \wid\ =\ \left [\ffi + \fsecdrv\right ]^{-1} 
                 \left [\wi^2 \wsecdrv + \wic\right ] 
\]
($\ffi\ \equal\ f(i\dt)$).
To prove convergence of the crystalline approximation
we shall derive the discrete analogue of estimates due to
Gage and Hamilton~\cite{Gage} for the curvature and then use standard
arguments to estimate $|\W(i\Delta\theta,t) - \wi(t)|$.
This will enable us to estimate the distance between 
the curve and the polygon.

        We mention that if the energy is isotropic, the initial
curve is a circle, and the initial polygon is regular and
circumscribed to the initial circle then the polygon remains
regular and circumscribed to the circle as they collapse to a point.
For other computational examples we refer the reader to Roosen
and Taylor~\cite{Roosen}.

        The organization of this \cpaper\ is as follows:
in \cSection~2 we set up the problem and give estimates for the weighted
curvature; then in \cSection~3 we prove the convergence of the 
crystalline algorithm for motion by weighted curvature.

        In this \cpaper\ we consider only convex curves in the plane.
A natural next step would be to consider general (i.e.\ nonconvex)
simple closed curves.
The paper by Grayson~\cite{Grayson} might be relevant
because it generalizes the results of~\cite{Gage} for
convex curves to general embedded plane curves.

\section{Setup and estimates for the weighted curvature}

Consider a smooth simple closed curve in the plane.  We say that
the curve is moving by weighted curvature when its normal velocity equals
its weighted curvature:
\be
        V(\theta,t)\ =\ \W(\theta,t)\ \equal\ 
                [f(\theta) + f^{\prime\prime}(\theta)]\ \kappa(\theta,t)
\label{vek}
\ee
as in Eq.~(4.11) of Angenent and Gurtin~\cite{AG}.
(This is Eq.~(41) of \cite{GK}.)  
Here $\theta$\/ is the angle between
the exterior normal and a fixed axis, 
$\kappa$\/ is the curvature, and
$\W$\/ is the weighted curvature.  We consider the case where $f$\/
(the interfacial energy per unit length) is positive, $C^3$, and such
that $f + \fpp\ >\ 0$.
We are interested in studying convergence of an approximation scheme
for this equation when the curve is convex.  In this scheme we
substitute the curve by a convex polygon with $N$\/ sides, $N\ \geq\ 4$.
The angle between two adjacent sides of the polygon is $\pi - \dt$,
where
\[
        \dt\ =\ \frac{\tpi}{N}
\]
is fixed, $0<\dt\leq\pi/2$.
The interior normals can be written 
$N_i\ \equal\ - (\cos i\dt, \sin i\dt)$
in a coordinate system that will henceforth remain fixed.  
We take $\theta$\/ to be the angle between the exterior 
normal and the $x$-axis.  Whenever we refer to ``sides''
we mean the sides of this polygon.
The $i$th side is the one with interior normal $N_i$\/ and
joins the $i$th vertex with the $(i+1)$th one.  We define
$T_i\ \equal\ (-\sin i\dt, \cos i\dt)$\/ (see Fig.~1).

\newpage

\includegraphics{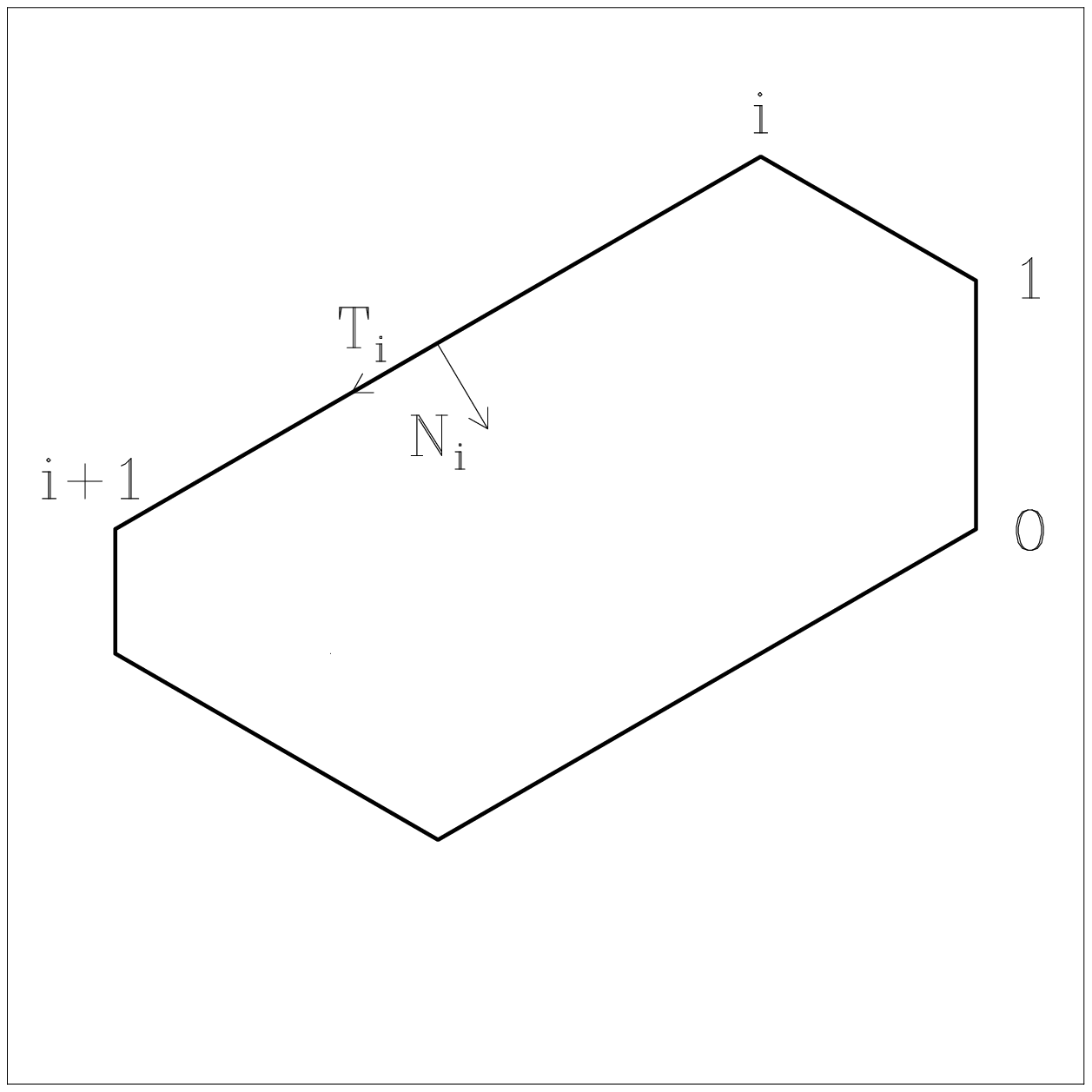}
\vglue 2.3truein                                              

\begin{center}
{\footnotesize\sc Fig. 1.}\label{ffive}  {\it \footnotesize
A polygon with six sides and labeled corners.} 
\end{center}

\vspace{2mm}

        Each side has zero curvature, of course, if 
curvature is defined in the standard way, namely the rate of change
with respect to arc length of the angle between the tangent 
and a fixed axis; in this standard sense the curvature is not defined at the
vertices of the polygon.  However, curvature can alternatively be defined as
the negative of the gradient of arclength;
weighted curvature is the negative of
the gradient of the interfacial energy.
In this sense each side of the polygon does have nonzero 
(weighted) curvature.
We denote by $L_i$, $V_i$, \ki, and \wi\ the length, velocity, curvature,
and weighted curvature, respectively, of the $i$th side.  
The normal velocity of the $i$th side in the approximation scheme is
\be
        V_i\ =\ \wi 
\label{vieki}
\ee
in the direction $N_i$.
Let $f_i\ \equal\ f(i\dt)$.  
As long as no sides disappear
the weighted curvature of the $i$th side is
\be
        \wi\ \equal\ \left [ f_i + 
                \frac{f_{i+1} - 2f_i + f_{i-1}}{\tomc} \right ]
                \ki
\label{ki}
\ee
and its curvature is
\be
        \ki\ \equal\ \fraction \frac{1}{\bigli}\ =\ 
                2 \tan\frac{\dt}{2}\,\frac{1}{\bigli}
\label{wi}
\ee
(see Eqs.~(45) of \cite{GK}).  
Also,
\be
        \dot{L}_i\ =\ \lidot
\label{lidot}
\ee
(see Eqs.~(10.18) of \cite{AG}).  

        Define
\[
\begin{array}{lclclcl}
        g&\equal&f +\ \fpp ,   &\qquad &\gi&\equal&\ffi + \fsecdrv,\\
         &      &              &       &   &      &                \\
        h&\equal&\d\frac{1}{g},&\m{and}&\hi&\equal&\d\frac{1}{\gi} .
\end{array}
\]
We note that $\gi\ >\ 0$\/ for all \dt, $0\ <\ \dt\ \leq\ \pi/2$,
and all $i$.  This is a consequence of the facts that (i)
$\gi\ \leq\ 0$\/ is equivalent to $\dt\ \neq\ \pi/2$\/ and
$\ffi\ \geq\ (\ffp + \ffm) / (2\cosdt)$, (ii) 
the line through $- N_{i-1} / \ffm$\/ and $- N_{i+1} / \ffp$\/
intersects the line with direction $- N_i$\/ at
$2\cosdt / (\ffp +  \ffm) \times (- N_i)$, and (iii) the curvature
of $(\cos(\ \cdot\ ),\sin(\ \cdot\ )) / f(\ \cdot\ )$\/ is
$f^3 (f + \fpp) / \sqrt{(f^2 + f^{\prime 2})^3}$\/ and hence
is positive, i.e.\ the polar diagram of $1 / f$\/ (the Frank
diagram) is convex.  Therefore $g_{\min}$\/ is bounded away
from zero; $g_{\max}$\/ is also bounded and these bounds are uniform
in \dt.  We use the notation $(\ \cdot\ )_{\min}$, $(\ \cdot\ )_{\max}$,
and $|\ \cdot\ |_{\max}$\/ 
for $\min_{0\leq i\leq N-1} (\ \cdot\ )_i$, 
$\max_{0\leq i\leq N-1} (\ \cdot\ )_i$,
and $\max_{0\leq i\leq N-1}\ |(\ \cdot\ )_i|$.

        Next we compute the rate of change of the length of the
polygon and of the area enclosed by it, and we compute the evolution
equation for the weighted curvature of the sides.  These computations
remain valid as long as no side disappears.  We could prove
directly from Eqs.~(\ref{ki}), (\ref{wi}), and (\ref{lidot}) that this
only happens when the polygon shrinks to a point or to a line
(see \cite{GK} 
or Taylor~\cite{TT}).
Instead we deduce it later from an estimate for the maximum
of the weighted curvature.

        The rate of change of the (total) length of the polygon is 
\begin{eqnarray}
\nonumber
        \dot{L}&=&\d \s \dot{L}_i\\
\nonumber&&\\
\nonumber
               &=&\d \s \left (\lidot\right )\\
\nonumber&&\\
\label{bigl}
               &=&\d \s 2 (\cotdt - \cscdt) \wi\\
\nonumber&&\\
\nonumber
               &=&\d - \s 4 \gi \frac{(\cotdt - \cscdt)^2}{\bigli}.
\end{eqnarray}
Throughout 
this \cpaper\ we write $\sum_i$\/ for $\sum_{i=0}^{N-1}$.
The last equation shows that, as in the continuous case, the total
length of the polygon is a decreasing function.  

        The area enclosed by the polygon can be computed by
\[
        A_{\Delta\theta}\ =\ \half \s d_i \bigli
\]
if we put the origin inside the polygon and let $d_i$ be the distance
between the origin and the line containing the $i$th side.
The rate of change of the area is
\[
\begin{array}{lcl}
        \dot{A}_{\Delta\theta}&=&\d \half \s \dot{d}_i \bigli +
                        \half \s d_i \dot{L}_i\\
&&\\
                &=&\d - \half \s \wi \bigli 
                        + \half \s d_i \left (\lidot\right )\\
&&\\
                &=&\d - \half \s \wi \bigli 
                   + \half \s \wi
                   \left (2 d_i \cotdt - d_{i+1} \cscdt - 
                                d_{i-1} \cscdt\right )\\
&&\\
                &=&\d - \s \wi \bigli\\ 
&&\\
                &=&\d - \s \gi \fraction 
                        \longrightarrow - \int_0^{2\pi} g(\theta)\, d\theta
\end{array}
\]
as $\Delta\theta\rightarrow 0$, since $\dot{d}_i= - \wi$\/ and   
$d_{i+1} \cscdt + d_{i-1} \cscdt - 2 d_i \cotdt = L_i$ 
(see Fig.~2).

\newpage

\includegraphics{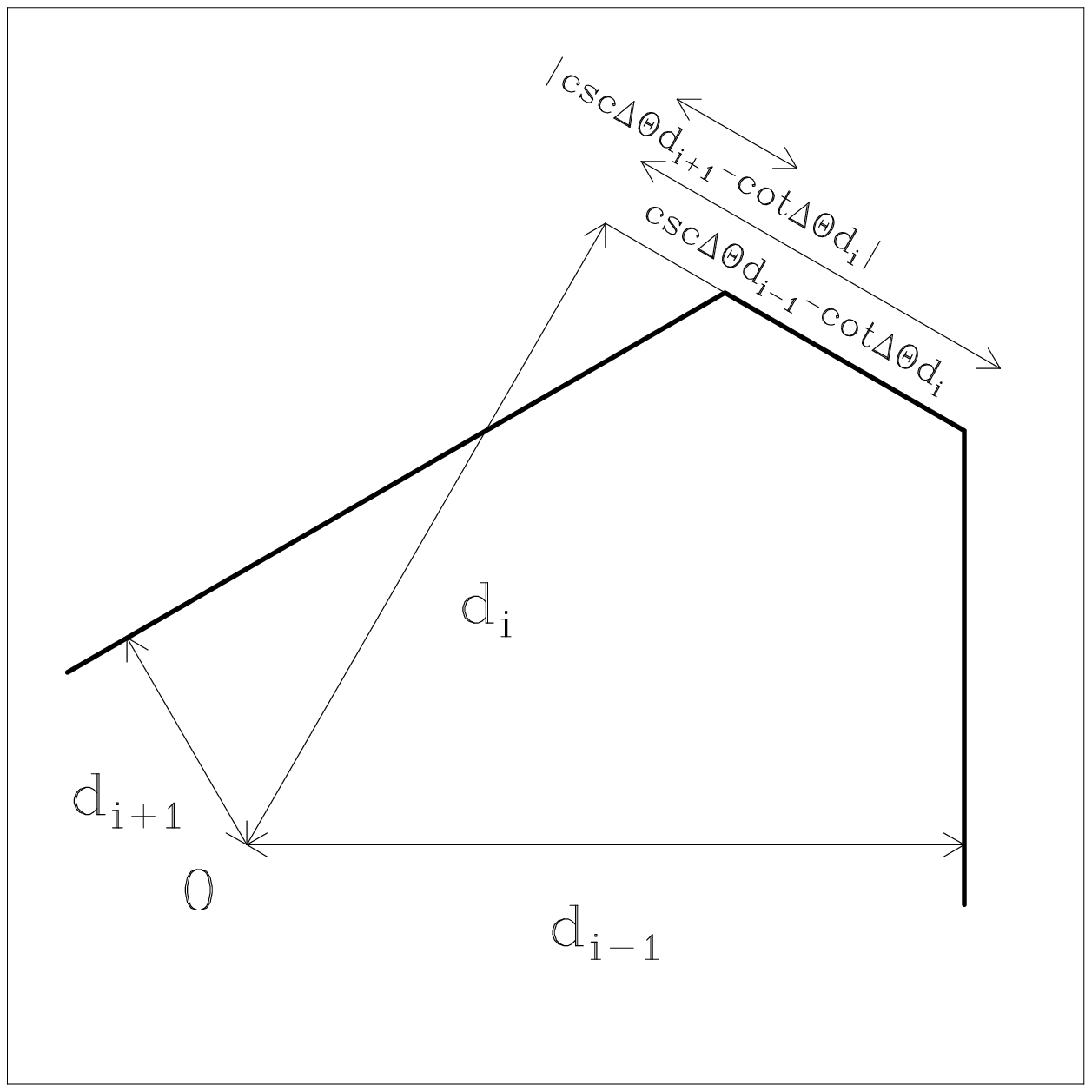}
\vglue 2.3truein                                              

\begin{center}
{\footnotesize\sc Fig. 2.}\label{fsix}  
{\it \footnotesize The $i$th and part of the $(i-1)$th and
        $(i+1)$th sides.  Here $N\ =\ 6$\/ and $i\ =\ 1$.}
\end{center}

\vspace{2mm}

By Eq.~(2.23) of Angenent and Gurtin~\cite{AG}
\be
        \W_t\ =\ h \left (\W^2 \W_{\theta\theta} + \W^3\right ) .
\label{kappat}
\ee
We compute the discrete analogue of Eq.~(\ref{kappat})
using Eqs.~(\ref{ki}), (\ref{wi}), and (\ref{lidot}):
\begin{eqnarray}
\nonumber
        \wid&=&\d - \gi \fraction \frac{\dot{L}_i}{{L_i}^2}\\
\nonumber&&\\
\nonumber
            &=&\d - \hi \ifraction \wis (\lidot)\\
\nonumber&&\\
\label{kidot}
            &=&\d \hi \left [\wis \wsecdrv + \wic\right ] .
\end{eqnarray}
This is a system of ordinary differential equations determining the
evolution of the $\wi$'s.  

        Note that one can reconstruct the polygon at time $t$
from the knowledge of the polygon at time zero and the $\wi$'s. 
In fact the (normal) velocity of the line containing the $i$th
side of the polygon is $\wi$.  Alternatively, it is easy to
check that the velocity of the midpoint of the $i$th side
is $\wi\, N_i\ -\ \frac{\omega_{i+1} - \omega_{i-1}}
{2\sin\Delta\theta}\, T_i$.

        Later we shall need the following discrete version of 
Poincar\'{e}'s inequality:
Let $\pp_0, \ldots, \pp_M$\/ be $M + 1$\/ real numbers with 
$\pp_0 = \pp_M = 0$.
Then
\be
        \sm \ppm^2\ \leq\ \frac{1}{2 (\omcpm)} \sm (\pmp - \ppm)^2 .
\label{poincare}
\ee
This is a restatement of the well know fact that the smallest
eigenvalue of the $(M-1)\times(M-1)$\/ matrix
\[\left [ \begin{array}{cccccc} 
        2       &-1     &0      &\cdots &0      &0\\ 
        -1      &2      &-1     &\cdots &0      &0\\ 
        \vdots  &\vdots &\vdots &\vdots &\vdots &\vdots\\ 
        0       &0      &0      &\cdots &2      &-1\\
        0       &0      &0      &\cdots &-1     &2
        \end{array} \right ]
\]
is $2 (\omcpm)$.
We remark that we will use the fact that the denominator
in Eq.~(\ref{kidot}) is \tomc\ and not $(\dt)^2$ because
inequality~(\ref{poincare}) is not valid if we substitute
$2(\omcpm)$ by $\frac{\pi^2}{M^2}$.

\vspace{\baselineskip}

        We now derive the discrete analogue of several estimates
by Gage and Hamilton \cite{Gage} for the curvature.  The computations
below are only valid as long as no side disappears.

\vspace{2mm}

\noindent {\bf First}, note that Eqs.~(\ref{kidot}) show that
$\w_{\min} (t)$ is a nondecreasing function.

\vspace{2mm}

\noindent {\bf Second}, we can bound the ``$H^1$ norm'' of 
the ``sequence \wi'' in terms of its ``$L^2$ norm.''  We compute
\[
\begin{array}{lcl}
        \ddt \s\!\! \left [ \wis - \wfdrvs \right ] \dt\!\!\!\!
         &=&\!\!\!\! 2 \d\s \left [ \wi \wid - \frac{(\wp - \wi)
                        (\dot{\w}_{i+1} - \dot{\w}_i)}
                {\tomc} \right ] \dt\!\!\!\\
&&\\
         &=&\!\!\!\! 2 \d\s \wid \left [ \wi + \wsecdrv \right ] \dt\\
&&\\
&=&\!\!\!\! 2 \d\s \hi \wis \left [ \wi + \wsecdrv \right ]^2\!\! \dt .
\!\!\!
\end{array}
\]
This is nonnegative so
\be
        \s \wfdrvs \dt\ \leq\ \s \wis \dt + c_1 .
\label{cone}
\ee
Here $c_1$\/ does not depend on time but only on the polygon at
time zero.  
Note that if $\w_{\max}(0)$\/ and $(\wp(0) - \wi(0))^2 / [\tomc]$\/
are uniformly bounded then $c_1$\/ is also uniformly bounded.
We may suppose $c_1\ \geq\ 0$.

\vspace{2mm}

\noindent {\bf Third}, if $A_{\Delta\theta}\ \geq\ \epsilon\ >\ 0$\/ then
we can bound 
\[
        \ws\ \equal\ 
\left\{\begin{array}{ll}
        \d\max_{0\leq j\leq N-1}\ \ \ \min_{j+1\leq i\leq j+N/2} \wi
&\qquad\m{for $N$\/ even,}
\\ &\\
        \d\max_{0\leq j\leq N-1}\ \ \ \min_{j+1\leq i\leq j+(N-1)/2} \wi
&\qquad\m{for $N$\/ odd,}
\end{array}\right . 
\]
which Gage and Hamilton call the ``median weighted curvature.''
Let $j_0$ be a value of $j$ for which the maximum is assumed.
A polygon with median weighted curvature \ws\ lies between parallel
lines whose distance is less than
\[
\begin{array}{lcl}
        \d \sum_{j=j_0}^{j_0+N/2} \sin((j-j_0)\dt) L_j
          &=&\d \sum_{j=1}^{N/2}\ \sin(j\dt) 
                g_{j_0+j} \fraction \frac{1}{\w_{j_0+j}}\\
&&\\
          &\leq&\d \frac{\sindt}{1 - \cosdt} g_{\max} \fraction \frac{1}{\ws}\\
&&\\
          &=&\d 2 \frac{g_{\max}}{\ws} 
\end{array}
\]
for $N$\/ even,
\[
\begin{array}{lcl}
        \d \sum_{j=j_0}^{j_0+(N-1)/2}\! \sin((j-j_0)\dt) L_j
        &=&\d \sum_{j=1}^{(N-1)/2}\ \sin(j\dt) 
                g_{j_0+j} \fraction \frac{1}{\w_{j_0+j}}\\
&&\\
        &\leq&\d 
        \frac{\sin\frac{\Delta\theta}{2}(\cos\frac{\Delta\theta}{2} + 1)}
             {1 - \cosdt} g_{\max} \fraction\frac{1}{\ws}\\
&&\\
        &=&\d (1 + \sec\frac{\Delta\theta}{2})
                \frac{g_{\max}}{\ws}\\
&&\\
        &\leq&\d (1 + \sec\frac{2\pi}{10})
                \frac{g_{\max}}{\ws}\ \leq\ \frac{5}{2} \frac{g_{\max}}{\ws} 
\end{array}
\]
for $N$\/ odd, and has diameter bounded by $L/2$.  So
\[
        A_{\Delta\theta}\ \leq\ \frac{5}{2} \frac{g_{\max}}{\ws} \frac{L}{2}\ 
                =\ \frac{5}{4} g_{\max} \frac{L}{\ws},
\]
or
\[
        \ws\ \leq\ \frac{5}{4} g_{\max} \frac{L}{A_{\Delta\theta}} .
\]
Note that if $A_{\Delta\theta}$\/ is bounded away from zero
uniformly in $\dt$\/ then $\ws$\/ is also uniformly bounded.

\vspace{2mm}

\noindent {\bf Fourth}, if \ws\ is bounded then
$\sum_i \gi \l \wi \dt$\/ is bounded.
By Eqs.~(\ref{kidot}),
\[
\begin{array}{lcl}
        \ddt \s \gi \l \wi&=&\s \left [ \wi \wsecdrv + \wis \right ]\\
&&\\
                      &=&\s \left [ \wis - \wfdrvs \right ].
\end{array}
\]
Define $\tilde{I}\ \equal\ \{ i\in \m{\bf N} | \wi\ >\ \ws \}$,
divide $\tilde{I}$\/ in maximal subsets of the form 
$\tilde{I}_j\ =\ \{\ijnew, \ijnew + 1, \ldots, \ijnew + \mj -2\}$, let
$I_j\ \equal\ \{\ijnew - 1, \ijnew, \ijnew + 1, \ldots, \ijnew + \mj -2\}$,
$I\ \equal\ \cup_j I_j$, and
$L\ \equal\ \m{\bf N}\setminus I$.  
Note that $I_j$\/ has \mj\ elements and note also that
$\mj - 1\ \leq\ N/2 -1$, or $\mj\ \leq\ N/2$,
because by the definition of \ws\ there are at most $N/2 -1$\/ \wi's
for $N$\/ even, $(N-1)/2 - 1$\/ \wi's for $N$\/ odd, 
corresponding to adjacent sides and with $\wi > \ws$.   
We have
\[
        \sum_{i\in L\cap\{0, \ldots, N-1\}}
                \left [ \wis - \wfdrvs \right ]
        \leq\ \sum_{i\in L\cap\{0, \ldots, N-1\}} \ws^2 
\]
and
\[
\begin{array}{lcl}
        \d\sum_{i\in I_j}\ \left [ \wis - \wfdrvs \right ]
         &\leq&\d \left [ \ws^2 - \frac{(\w_{i_j} - \ws)^2}{\tomc} \right ]\\
&&\\
         &    &\d + \sum_{i=i_j}^{i_j+M_j-3}\ \left [ \wis - \wfdrvs \right ]\\
&&\\
         &    &\d + \left [ \w_{i_j+M_j-2}^2 - 
                \frac{(\ws - \w_{i_j+M_j-2})^2}{\tomc} \right ]\\
&&\\
         &\leq&\d 2 \ws \sum_{i=i_j}^{i_j+M_j-2}\ \wi.
\end{array}
\]
For the last inequality we have used Eq.~(\ref{poincare}) with $M\ =\ \mj$,
$\ppm\ =\ \w_{i_j-1+m} - \ws$\/ for $1 \leq m \leq \mj - 1$, the fact that
\[
        \d\frac{1}{2 (1 - \cos \frac{\pi}{M_j} )}\ \leq\ 
        \d\frac{1}{2 (1 - \cos \frac{2\pi}{N} )}\ =\
        \d\frac{1}{\tomc},
\]
and also that $\mj \geq 2$.  The last estimates together with Eq.~(\ref{bigl})
imply
\[
        \ddt \s \gi \l \wi \dt\ \leq\ \ws^2 \frac{\tpi}{\dt} \dt - 
                2 \ws \ifraction \dt \s \dot{L}_i .
\]
Thus if $\ws(\tau) \leq \Omega$\/ for $0\ \leq\ \tau\ \leq\ t$\/ then
\[
        \s \gi \l \wi(t) \dt\ \leq\ \s \gi \l \wi(0) \dt +
                \tpi {\Omega}^2 t + 2 \Omega [L(0) - L(t)] .
\]
Note that if $\w_{\max}(0)$\/ and $\Omega$\/ are 
uniformly bounded then this bound is also uniform.

\vspace{2mm}

\noindent {\bf Fifth}, if $\sum_i \gi \l \wi \dt$\/ is bounded then for any 
$\delta>0$\/ we can find a constant $c_2\ >\ 1$\/ such that $\wi \leq c_2$\/
except for $i\dt$\/ in intervals whose total length is less than $\delta$.
In fact, $\wi\geq c_2$\/ for $L$ values of $i$ and $L\dt\geq\delta$\/ 
implies
\be
\begin{array}{lcl}
        \s \gi \l \wi \dt&\geq&
                g_{\min}\dt L\,\l c_2 + g_{\max} (N - L) \dt\,\l \w_{\min}\\
        &\geq&g_{\min}\delta\,\l c_2 + g_{\max}(\tpi - \delta)\,\l \w_{\min} ,
\end{array}
\label{kmin}
\ee
and this gives a contradiction when $c_2$ is large (we have assumed
$\w_{\min}\ <\ 1$).  Here the
constant $c_2$ depends on $\w_{\min}(0)$.  
Note that if $\w_{\min}(0)$\/ is bounded away from zero 
uniformly in $\dt$\/ 
then the bound for $c_2$\/ is also uniform.

\vspace{2mm}

\noindent {\bf Sixth} and last, if $\wi \leq c_2$\/ except 
for $i\dt$\/ in intervals
whose total length is less than $\delta$\/ and $\delta$\/ is
small enough then $\w_{\max}$\/ is bounded.  
In fact, $\w_m \leq c_2$\/ and $n-m\leq\delta/\dt$\/ imply
\[
\begin{array}{lcl}
        \w_n&=&\d \w_m + \sum_{i=m}^{n-1}\ (\wp - \wi)\\
&&\\
           &\leq&\d c_2 + 
                \left [ \sum_{i=m}^{n-1}\ \frac{\tomc}{\dt} \right ]^{1/2}
                \left [ \sum_{i=m}^{n-1}\ \wfdrvs \dt \right ]^{1/2}\\
&&\\
           &\leq&\d c_2 + \sqrt{\delta} \closeone
                \left [ \sum_{i=m}^{n-1}\ \wfdrvs \dt \right ]^{1/2}\\
&&\\
           &\leq&\d c_2 + \sqrt{\delta}
                        \left [ \s \wis \dt + c_1 \right ]^{1/2} .
\end{array}
\]
Hence
\[
        \w_{\max}\ \leq\ c_2 + \sqrt{\delta} 
                \left [ \sqrt{\tpi} \w_{\max} +\sqrt{c_1} \right ] ,
\]
and for small $\delta$\/ we get $\w_{\max} \leq 2 c_2$.

        Consider a fixed polygon.  Let
\be
        T_{\Delta\theta}\ \equal\ A_{\Delta\theta}(0) \left / 
                \s \left \{ \gi \fraction\ \right \} \right . 
\label{time}
\ee
with $A_{\Delta\theta}(0)$\/ the area enclosed by the polygon at 
time zero.  Recall that the denominator is $-\dot{A}_{\Delta\theta}$.
Suppose that a side of the polygon disappears for 
$t\ <\ T_{\Delta\theta}$\/
and let $\bar{T}$\/ be the first time that happens.
If $0\ \leq\ t\ \leq\ \bar{T}$
then $A_{\Delta\theta}(t)$\/ is bounded away from zero.
Since the initial polygon is convex $\min_{0\leq i\leq N-1}\ \w_i(0)$\/ 
is positive.  The estimates above imply that 
$\sup_{0\leq t\leq \tilde{T}}\ 
\max_{0\leq i\leq N-1}\ \w_i(t)$\/ is bounded.
So $L_{\min}(\bar{T})\ >\ 0$.
This is a contradiction.  Therefore
the polygon will disappear in exactly $T_{\Delta\theta}$\/ 
units of time and no sides vanish before that.  
\setcounter{transfer}{\theequation}

\section{Convergence of the crystalline algorithm for motion by 
weighted curvature}

\setcounter{equation}{\thetransfer}
We wish to compare a convex polygon
moving according to Eqs.~(\ref{vieki}) with a smooth simple closed
convex curve moving according to Eq.~(\ref{vek}).  
To do this we
want to estimate the distance between the line containing the $i$th
side of the polygon and the tangent to the curve parallel to
this line.  Since the (normal) velocities of these lines 
at time $t$ are $\wi(t)$\/ and $\W_i(t)\ \equal\ \W(i\dt,t)$,
respectively, we estimate the difference $\wmw$.
All derivatives of $\W$\/ remain bounded for 
$0\ \leq\ t\ \leq\ \tilde{T}\ <\ T$\/ where
$T$\/ is the time that the area enclosed by the curve vanishes (see
Gage and Hamilton~\cite{Gage}).  So, for these values of $t$,
$\W_{\theta\theta\theta\theta}$\/ is bounded and
\[
        \W_{\theta\theta}(i\dt,t)\ =\ \Wsecdrv + \odts .
\]
Using this in Eq.~(\ref{kappat}) and subtracting Eqs.~(\ref{kidot}),
\[
\begin{array}{lcl}
        \dot{\W}_i - \wid
           &=&\d \hi \wis \frac{(\Wp - \wp) - 2 (\wmw) + (\Wm - \wm)}{\tomc}\\
&&\\
            & &\d + \hi (\wmw) (\Wi + \wi) \Wsecdrv\\
&&\\
            & &\d + \hi (\wmw) (\Wi^2 + \Wi\wi + \wis) + \odts.
\end{array}
\]
Let $\li\ \equal\ \wmw$.  The last equation can be written
\be
        \dot{\Lambda}_i\ =\ a_i \secdrvl + b_i \li + \odts
\label{appropriate}
\ee
with $a_i$\/  positive, and $b_i$\/ uniformly bounded
if $\w_{\max}$\/ is uniformly bounded.  Therefore
\[
        \dot{\Lambda}_{\max}\ \leq\ b \Lambda_{\max} + \odts
                \qquad\m{for}\ \ \Lambda_{\max}\ >\ 0 ,
\]
\[
        \dot{\Lambda}_{\min}\ \geq\ b \Lambda_{\min} + \odts
                \qquad\m{for}\ \ \Lambda_{\min}\ <\ 0 ,
\]
for $b = \sup_{0\leq t\leq \tilde{T}}\ \max_i |b_i|$.  So,
\be
        |\Lambda|_{\max}(t)\ \leq\ |\Lambda|_{\max}(0) e^{bt}
                        + \odts\, t e^{bt} 
\label{lmax}
\ee

        Two important tasks remain to be done.  One is to specify how we
obtain the initial polygon and derive estimates for the
curvature at time zero.  The other one is to prove that if $\wmw$\/
is small then the curves are in fact uniformly close to one another.
We address these issues in sequence.
\addtocounter{equation}{+1}
\setcounter{initialplyg}{\theequation}

\vspace{2mm}
\noindent $(\theinitialplyg_{\Delta\theta})$

\vspace{1mm}

{\em The initial polygon is obtained by the union of segments
on lines tangent 

to the initial curve at
the points with exterior normal}\/
$(\cos i\dt, \sin i\dt)$.

\vspace{2mm}

\noindent It should be emphasised that there should be a segment
of nonzero length corresponding to each $i$.
We derive estimates for $|\Lambda|_{\max}(0)$\/ and $|\Upsilon|_{\max}(0)$,
\[
        \upi(0)\ \equal\ \W_{\theta}(i\dt,0) - \frac{\wp(0) - \wi(0)}{\sindt}.
\]
We need an estimate for $|\Upsilon|_{\max}$\/ to prove that
$c_2$\/ (the constant in Eq.~(\ref{kmin})) remains bounded 
as $\dt \rightarrow 0$.
We prove that
$\ki(0)\ =\ \kpi(0) + \odts$\/ 
and $(\kp(0) - \ki(0)) / \sindt\ =\ \kappa_{\theta}(i\dt,0) + \odt$\/
($\kpi(t)\ \equal\ \kappa(i\dt,t)$).  Since
\be
        \li(0)\ =\ \Wi(0) - \wi(0)\ =\ \gi (\kpi(0) - \ki(0))
\label{lambdai}
\ee
and
\be
\begin{array}{lcl}
        \upi(0)&=&\d\gi \kappa_{\theta}(i\dt,0) 
                        + g_{\theta}(i\dt) \kpi(0)\\
&&\\
               & &\d - \gp \frac{\kp(0) - \ki(0)}{\sindt}
                        - \frac{\gp - \gi}{\sindt} \ki(0) ,
\end{array} 
\label{upsiloni}
\ee
it will follow that $|\Lambda|_{\max}(0)\ =\ \odts$\/ 
and $|\Upsilon|_{\max}(0)\ =\ \odt$.
To estimate $\ki(0)$\/ we estimate the length $L_i(0)$.
We decompose $L_i(0)$\/ into $L_i(0) = \lip + \limn$\/
(see Fig.~3).

\includegraphics{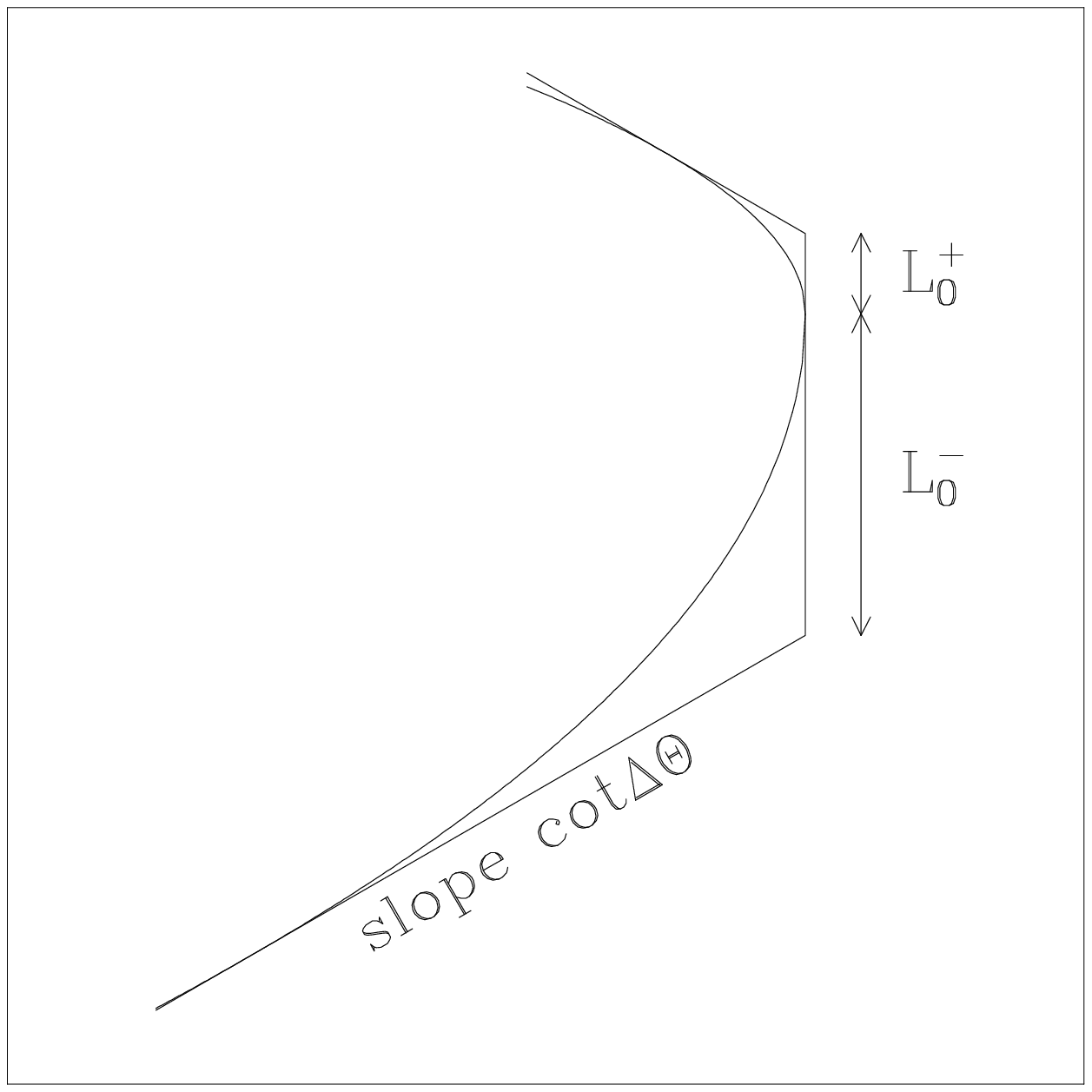}
\vglue 2.3truein                                              

{\footnotesize\sc Fig. 3.}\label{feight}  
{\it \footnotesize The $i$th and part of the $(i-1)$th and
$(i+1)$th sides at time zero.  $L_i(0) = \lip + \limn$.
Here $N\ =\ 6$\/ and $i\ =\ 0$.}

\vspace{2mm}

Note that
\[
        \kappa(\theta,0)\ =\ \kpi(0) + \tmidt \kptidt + \otmidts
\]
and
\[
\begin{array}{lcl}
        r(\theta,0)\!\!&=&\!\d r(i\dt,0) + \int_{i\Delta\theta}^{\theta} 
                        \frac{T(\vartheta)}{\kappa(\vartheta,0)} \dvt\\
&&\\
                 &=&\!\d r(i\dt,0) + \int_{i\Delta\theta}^{\theta} 
                        \frac{(-\sin\vartheta, \cos\vartheta)}
                    {\kpi(0) + (\vartheta - i\dt) \kptidt + \ovtmidts} \dvt .
\end{array}
\]
Therefore
\[
        \lip\ =\ \ikpi \firstint - \cotdt \ikpi \secint .
\]
where $\rho_i\ =\ \kappa_{\theta}(i\dt,0) / \kappa_i(0)$.
One easily checks that
\[
        \cotdt\ =\ \frac{1}{\dt} - \frac{\dt}{3} + \odtc ,
\]
\be
\begin{array}{rcl}
        \firstint&=&\d \dt - \frac{(\dt)^2}{2} \rho_i + \odtc ,\\
&&\\
        \secint&=&\d \frac{(\dt)^2}{2} - \frac{(\dt)^3}{3} \rho_i + \odtf .
\end{array}
\label{middle}
\ee
So
\be
        \lip\ =\ \ikpi \left [ \frac{\dt}{2} - \frac{(\dt)^2}{6} \rho_i
                + \odtc \right ] .
\label{length}
\ee
But then, of course,
\[
        \limn\ =\ \ikpi \left [ \frac{\dt}{2} + \frac{(\dt)^2}{6} \rho_i
                + \odtc \right ] .
\]
Adding, the expression for $L_i(0)$\/ in terms of \kpi(0) and \dt\ is
\[
        L_i(0)\ =\ \ikpi \left [ \dt + \odtc \right ] .
\]
This leads to the desired expression of $\ki(0)$\/ in terms of $\kpi(0)$,
\begin{eqnarray}
\nonumber
        \ki(0)&=&\d \fraction \frac{1}{L_i(0)}\\
\nonumber&&\\
\nonumber
           &=&\d \frac{\dt + \odtc}{\frac{1}{\kappa_i(0)} \left [
                        \dt + \odtc \right ]}\\
\nonumber&&\\
\label{errork}
           &=&\d \kpi(0) + \odts ,
\end{eqnarray}
and $(\kp(0) - \ki(0)) / \sindt$\/ in terms of $\kptidt$,
\begin{eqnarray}
\nonumber
        \d \frac{\kp(0) - \ki(0)}{\sindt}
           &=&\d \frac{\left \{\kpp(0) + \odts\right \} -
                       \left \{\kpi(0) + \odts\right \}}
                      {\dt + \odtc}\\
\nonumber&&\\
\nonumber
           &=&\d \frac{\kpp(0) - \kpi(0)}{\dt} + \odt\\
\nonumber&&\\
           &=&\d \kappa_\theta (i\dt,0) + \odt .
\label{errorv}
\end{eqnarray}  

        Now we prove that the curves are uniformly close to one another.
It is convenient to use the Hausdorff metric on compact sets $A$\/ and $B$,
defined as
\[
        D(A,B)\ \equal\ \max    \left \{
                        \max_{p\in A} \min_{q\in B} \m{dist}(p,q),
                        \max_{q\in B} \min_{p\in A} \m{dist}(p,q)
                                \right \}.
\]
Thus, $A$\/ and $B$\/ are close if any point of $A$\/ is
close to $B$\/ and vice-versa.
Our main result is the following

{\sc Theorem.}  
{\it Suppose $f$\/ $($the interfacial energy$)$
is positive, $C^3$, and $f + \fpp\ >\ 0$.
Consider a smooth simple closed convex
curve $C(t)$\/ enclosing at time zero an area $A$\/ and 
moving according to\/ {\rm Eq.~$($\ref{vek}$)$}.  Consider further a polygon $P(t)$\/ 
determined initially by condition~$(\theinitialplyg_{\Delta\theta})$\/ and
moving according to\/ {\rm Eqs.~$($\ref{vieki}$)$}.  
Let $T\ \equal\ A\ /
\int_0^{2\pi} (f(\theta) + \fpp(\theta))\, d\theta$\/
and\/ choose any $\tilde{T}\ <\ T$\/
$($\/$T$\/ is the time it takes for the curve to shrink to a point or 
to a straight line$)$.
If\/ $0\ \leq\ t\ \leq\ \tilde{T}$\/ 
and $\dt$\/ is sufficiently small then $P(t)$\/ is close to 
$C(t)$\/ in the sense that} $D(P(t),C(t))\ =\ \odts${\it .}
\begin{proof}
By Eqs.~(\ref{lambdai}), (\ref{upsiloni}), 
(\ref{errork}), and (\ref{errorv}),
condition~$(\theinitialplyg_{\Delta\theta})$\/ implies that
\[
        |\Lambda|_{\max}(0)\ =\ \max_{0\leq i\leq N-1}\ 
                |\W(i\dt, 0) - \wi(0)|\ =\ \odts
\]
and
\[
        |\Upsilon|_{\max}(0)\ =\ \max_{0\leq i\leq N-1}\ 
        \left |\W_{\theta}(i\dt, 0) - \frac{\wp(0) - \wi(0)}
                                                        {\sindt}\right |
\ =\ \odt.
\]
These estimates imply that $\w_{\min}(0)$\/ is bounded away from
zero, $\w_{\max}(0)$\/ is bounded, and $(\wp(0) - \wi(0))^2 / [\tomc]$\/
is bounded with these bounds uniform in $\dt$.  
Pick $\dt$\/ sufficiently small so that
$T_{\Delta\theta}\ >\ \tilde{T} + \delta$\/ for some positive $\delta$\/
and $A_{\Delta\theta}(\tilde{T})$\/ is bounded away from zero uniformly
in $\dt$\/ ($T_{\Delta\theta}$\/ is defined in Eq.~(\ref{time})).  
Then the bounds of \cSectiontwo\ are uniform in $\dt$.
Since $0\ \leq t\ \leq\ \tilde{T}\ <\ T$, 
$\W_{\theta\theta\theta\theta}$\/
is uniformly bounded.  Then Eq.~(\ref{lmax}) implies
\be
        |\W(i\dt, t) - \wi(t)|\ \leq\ c_3 \dts .
\label{square}
\ee

        We define $C_P(t)$\/ to be the polygon obtained by the union
of segments on lines tangent to $C(t)$\/ at points with exterior
normal $-N_i$.  We prove that $D(P,C)\ =\ \odts$\/ in two steps.
First we estimate $D(P,C_P)$\/ and then we estimate $D(C,C_P)$.

        We claim that the distance $D(P,C_P)\ =\ \odts$.
The (normal) velocity of the line containing the $i$th
side of the polygon $P$\/ is $\wi$\/ while the (normal) velocity of
the tangent to the curve parallel to this line is $\W_i$.
Since these lines coincide at time zero inequality~(\ref{square})
implies that their distance at time $t$\/ is less or equal to
$c_3 t \dts$.  Furthermore, this is true for all $i$\/ so
$P$\/ lies in a strip of width $2 c_3 t \dts$\/ whose center is $C_P$.
The Hausdorff distance between $C_P$\/ and the boundary of this
strip is $c_3 t \dts / \cos(\dt/2)$.  Hence
the distance $D(P,C_P)\ =\ \odts$.    

        Now we claim that the distance $D(C,C_P)\ =\ \odts$.
This distance is equal to the distance between 
some vertex of $C_P$\/ and the point on $C$\/ closest to it,
because $C$\/ is convex.  If the vertex in question is the $(i+1)$th
one, then $D(C,C_P)$ is less than the distance between the
$(i+1)$th vertex of $C_P$\/ and the point on $C$\/ with exterior normal
$(\cos(i+1/2)\dt, \sin(i+1/2)\dt)$.
The point on $C$\/ with exterior normal
$(\cos(i+1/2)\dt, \sin(i+1/2)\dt)$\/ is
\[
        p_i + \int_{i\Delta\theta}^{(i+1/2)\Delta\theta} \frac{T(\vartheta)}
                {\kappa(\vartheta,t)} \dvt ,
\]
where $p_i$\/ is the point of tangency between $C$\/ and the
side of $C_P$\/ with exterior normal $-N_i$, or
\[
        p_i + \firstintm\ T_i + \secintm\ N_i ;
\]
by Eqs.~(\ref{middle}) this equals
\[
        p_i + \left [ \frac{\dt}{2\kpi} + \odts \right ] T_i + \odts\, N_i .
\]
On the other hand, by Eqs.~(\ref{length}) the $(i+1)$th vertex of $C_P$\/ is
\[
        p_i +  \left [ \frac{\dt}{2\kpi} + \odts \right ] T_i .
\]
We conclude that the distance between the
$(i+1)$th vertex of $C_P$\/ and the point on $C$\/ with normal
$(\cos(i+1/2)\dt, \sin(i+1/2)\dt)$\/ is \odts.

        Finally, 
$$
        D(P,C)\ \leq\ D(P,C_P) + D(C,C_P)\ =\ \odts . 
$$
\end{proof}

\vspace{\baselineskip}

        We close with two remarks.  The first one concerns the
comparison between $\W$\/ and its discrete analogue.  Above we
proved a result of convergence in the sup norm. 
It is also possible to prove that 
\[
        \left |\W_{\theta}(i\dt, t) - \frac{\wp(t) - \wi(t)}
                                        {\sindt}\right |\ =\ \odt ,
\]
for\/ $0\ \leq\ t\ \leq\ \tilde{T}$.
The proof goes in three steps: (i) one derives the evolution equation
for $(\wp - \wi) / \sindt$\/ using Eq.~(\ref{kidot}), (ii)
one differentiates Eq.~(\ref{kappat}) with respect to $\theta$\/
and writes an analogue of the equation obtained in the previous step
by substituting derivatives of $\W_{\theta}$\/ by difference quotients,
and (iii) one subtracts the equations obtained in the previous
two steps, rewrites the difference in an appropriate way (as in
Eq.~(\ref{appropriate})), and derives a differential inequality
for the maximum of the absolute value of the difference we want to
estimate.  One estimates $(\wp - \wi) / \sindt$\/ 
and not $(\wp - \wm) / (2\sindt)$\/ directly to be able to
make conclusions about the sign of the difference quotients that appear
in step three.

        Our final remark concerns the comparison between the order of
convergence of this numerical scheme for simple closed convex curves,
on the one hand, and for graphs with prescribed boundary conditions,
on the other.  The main result of 
this \cpaper\  
is that the former order of convergence is
$\odts$\/ in the Hausdorff norm, when the admissible angles are 
equally spaced. The latter one is
$\mbox{O}(\sqrt{\Delta\theta})$\/ in the stronger norm $H^1$,
for a general distribution of admissible slopes
(see \cite{GK}).

\vspace{3mm}

{\bf Acknowledgment.} 
This paper is part of the author's Ph.D. thesis.
The author wishes to thank Robert V. Kohn for his advice.

\end{document}